\documentclass[14pt]{article}
\usepackage{mathrsfs}
\usepackage{amsthm}
\usepackage{amsmath}
\usepackage{graphicx}
\usepackage{color}
\usepackage{amsfonts}
\newtheorem{theorem}{Theorem}[section]

\newtheorem{lemma}[theorem]{Lemma}

\numberwithin{equation}{section}
\normalsize

\begin{document}
\title{\textbf{Contact processes with random recovery rates and edge weights on complete graphs}}

\author{Xiaofeng Xue \thanks{\textbf{E-mail}: xfxue@bjtu.edu.cn \textbf{Address}: School of Science, Beijing Jiaotong University, Beijing 100044, China.}\\ Beijing Jiaotong University \\ Yu Pan \thanks{\textbf{E-mail}: perryfly@pku.edu.cn \textbf{Address}: School of Mathematical Sciences, Peking University, Beijing 100049, China.}\\Peking University}

\date{}
\maketitle

\noindent {\bf Abstract:}
In this paper we are concerned with the contact process with random recovery rates and edge weights on complete graph with $n$ vertices. We show that the model has a
critical value which is inversely proportional to the product of the mean of the edge weight and the mean of the inverse of the recovery rate.
In the subcritical case, the process dies out before a moment with order $O(\log n)$ with high probability as $n\rightarrow+\infty$. In the supercritical case,
the process survives at a moment with order $\exp\{O(n)\}$ with high probability as $n\rightarrow+\infty$. Our proof for the subcritical case is inspired
by the graphical method introduce in \cite{Har1978}. Our proof for the supercritical case is inspired by approach introduced in \cite{Pet2011}, which deal with the case
where the contact process is with random vertex weights.

\quad


\section{Introduction}\label{section 1}
In this paper we are concerned with contact processes with random recovery rates and edge weights on complete graphs. For each integer $n\geq 1$, we denote by $C_n$ the complete graph with $n$ vertices. We denote by $1,2,3,\ldots,n$ the $n$ vertices of $C_n$ and $(i,j)$ the edge connecting $i$ and $j$ for $1\leq i\neq j\leq n$. Hence, $(i,j)=(j,i)$ for each pair of $i \neq j$. Under our notations, $C_n$ can be seen as a subgraph of $C_m$ when $n<m$.

Let $\rho$ and $\xi$ be two random variables such that $P(0\leq\rho\leq 1)=1, P(\rho>0)>0$ and $P(1\leq\xi\leq M)=1$ for some $M\in (1,+\infty)$, then we assume that $\{\xi(i)\}_{i=1}^{+\infty}$ are i. i. d. copies of $\xi$ while $\{\rho(i,j):1\leq i<j\}$ are i. i. d. copies of $\rho$ and independent of $\{\xi(i)\}_{i=1}^{+\infty}$. For $i>j$, we define $\rho(i,j)=\rho(j,i)$. The contact process on $C_n$ with random recovery rates $\{\xi(i)\}_{1\leq i\leq n}$ and edge weights $\{\rho(i,j)\}_{1\leq i\neq j\leq n}$ is a spin system with state space $X=\{0,1\}^{\{1,2,\ldots,n\}}$ and generator function $\mathcal{L}$ given by
\begin{equation}\label{equ 1.1 generator}
\mathcal{L}f(\eta)=\sum_{i=1}^nc(\eta,i)\big(f(\eta^i)-f(\eta)\big)
\end{equation}
for any $\eta\in X$ and $f\in C(X)$, where
\[
\eta^i(j)=
\begin{cases}
\eta(j) & \text{~if~} j\neq i,\\
1-\eta(i) & \text{~if~} j=i
\end{cases}
\]
and
\begin{equation}\label{equ 1.2 spin function}
c(\eta,i)=
\begin{cases}
\xi(i) &\text{~if~}\eta(i)=1,\\
\frac{\lambda}{n}\sum\limits_{1\leq j\leq n,\atop j\neq i}\rho(i,j)\eta(j) &\text{~if~} \eta(i)=0,
\end{cases}
\end{equation}
where $\lambda$ is a positive parameter called the infection rate. For any $t\geq 0$, we denote by $\eta_t$ the configuration of the process at moment $t$.

Intuitively, the process describes the spread of an infectious disease on $C_n$. $\{i:\eta(i)=1\}$ are infected vertices while $\{i:\eta(i)=0\}$ are healthy vertices. When vertex $i$ is infected, it waits for an exponential time with rate $\xi(i)$ to become healthy. When $i$ is healthy while $j$ is infected, then $j$ infects $i$ at rate proportional to $\rho(i,j)$, which is the weight on the edge connecting $i$ and $j$.

When $\xi=\rho=1$, then our model is the classic contact process introduced by Harris in \cite{Har1974}. Please see Chapter 6 of \cite{Lig1985} and Part one of \cite{Lig1999} for a detailed survey of the classic contact process. When $\rho=1$, then our model is the contact process with random recovery rates, the study of which dates back to 1990s. In \cite{Bram1991}, Bramson, Durrett and Schonmann study the contact process on $\mathbb{Z}^1$ where each vertex is in a bad situation with probability $p$. A vertex in bad situation has recovery rate smaller than that in good situation. They show that the model has a intermediate phase in which the process survives but does not grow linearly. In \cite{Lig1992}, Liggett studies contact process on $\mathbb{Z}^1$ with general random recovery rates and gives a sufficient condition for the process to survive. When $\xi=1$, our model is the contact process with random edge weights. Especially, when $\xi=1$ and $P(\rho=1)=p=1-P(\rho=0)$, then our model is contact process on clusters of bond percolation. In \cite{Chen2009}, Chen and Yao prove that complete convergence theorem holds for contact process on clusters of bond percolation on $\mathbb{Z}^d\times \mathbb{Z}^+$. In \cite{Yao2012}, they extend this result to the case where the process with general random edge weights on $\mathbb{Z}^+\times \mathbb{Z}^d$. In \cite{Xue2016}, Xue shows that contact process on clusters of oriented bond percolation on $\mathbb{Z}^d$ has critical value $(1+o(1))/(dp)$ for large $d$, where $p$ is probability that an given edge is open. The conclusion in \cite{Xue2016} is easy to extend to the case where the process is with general random edge weights on oriented lattice. The bond percolation on complete graph is also known as Erdos-Renyi model (see Chapter 3 of \cite{van2012}). Our model contains contact process on Erdos-Renyi graph $G(n,p)$ as a special case when $P(\rho=1)=p=1-P(\rho=0)$ and $\xi=1$.

In this paper we assign i. i. d. infection weights on edges. It is also interesting to assign weights on vertices and let an infectious vertex infect a healthy one at rate proportional to the production of the vertex-weights of the two vertices. In \cite{Pet2011}, Peterson studies such contact process with random vertex weights on complete graphs $C_n$. It is shown in \cite{Pet2011} that the model has critical value $\lambda_c$  which is the inverse of the second moment of the vertex weight. When $\lambda>\lambda_c$, then the process survives at moment $e^{Cn}$ for some $C>0$ with high probability. When $\lambda<\lambda_c$, then the process dies out before moment $c\log n$ for some $c>0$ with high probability. We are inspired by \cite{Pet2011} a lot. The main motivation of this paper is proving counterpart of the above main result in \cite{Pet2011} for the case where the process is with random edge weights and recovery rates. Inspired by \cite{Pet2011}, Xue studies contact process with random vertex weights on oriented lattice $\mathbb{Z}^d$ in \cite{Xue2015} and shows that the process has critical value $(1+o(1))/(d{\rm E}\rho^2)$ for large $d$, where $\rho$ is the vertex weight. When the vertex weight belongs to $\{0,1\}$, then contact process with random vertex weights turns into the process on clusters of site-percolation, which is studied in \cite{Ber2011}, where several sufficient conditions for the process to survive are given.

This paper is arranged as following. In Section \ref{section 2}, we give our main result, which is the counterpart of the main result in \cite{Pet2011} for the process with random recovery rates and edge weights. In Section \ref{section 3}, we give an intuitive explanation of our main result according to a mean-field analysis. In Section \ref{section 4}, we give the proof of the subcritical case. Our proof is inspired by the approach of graphical representation introduced in  \cite{Har1978} by Harris. In Section \ref{section 5}, we give the proof of the supercritical case. Our proof is inspired by that in \cite{Pet2011}.

\section{Main results}\label{section 2}

In this section we give our main result. First we introduce some notations. We assume that $\{\xi(i)\}_{i=1}^{+\infty}$ and $\{\rho(i,j)\}_{1\leq i\neq j<+\infty}$ are defined under probability space $\{\Omega,\mathfrak{F},\mu\}$. We denote by ${\rm E}_\mu$ the expectation operator with respect to $\mu$. For any sample $\omega\in \Omega$, we denote by $P_{\lambda,n}^{\omega}$ the probability measure of the contact process on $C_n$ with random recovery rates $\{\xi(\omega,i)\}_{i=1}^n$, edge weights $\{\rho(\omega,i,j)\}_{1\leq i\neq j\leq n}$ and infection rate $\lambda$. $P_{\lambda,n}^{\omega}$ is called the quenched measure. We denote by ${\rm E}_{\lambda,n}^{\omega}$ the expectation operator with respect to $P_{\lambda,n}^{\omega}$. We define
\[
P_{\lambda,n}(\cdot)=\int P_{\lambda,n}^\omega(\cdot)~\mu(d\omega)={\rm E}_\mu\big(P_{\lambda,n}^\omega(\cdot)\big),
\]
which is called the annealed measure. We denote by ${\rm E}_{\lambda,n}$ the expectation operator with respect to $P_{\lambda,n}$. When $\{i:\eta_0(i)=1\}=A$, then we write $\eta_t$ as $\eta_t^A$, but we omit the superscript when $A=C_n$. For simplicity, we identify $\eta_t$ with the set $\{i:\eta_t(i)=1\}$ when there is no misunderstanding.

Now we can give our main result. Assuming that all the vertices are infected at $t=0$, then our model performances the following phase transition phenomenon.

\begin{theorem}\label{theorem 2.1 main}
When $\lambda<\lambda_c=\frac{1}{{\rm E}\rho{\rm E}\frac{1}{\xi}}$, then there exists $c(\lambda)>0$ such that
\begin{equation}\label{equ 2.1 subcritical}
\lim_{n\rightarrow+\infty}P_{\lambda,n}^{\omega}(\eta_{c(\lambda)\log n}=\emptyset)=1 \text{\quad a.s.}
\end{equation}
with respect to the probability measure $\mu$. When $\lambda>\lambda_c$, then there exists $C(\lambda)>0$ such that
\begin{equation}\label{equ 2.2 supercritical}
\lim_{n\rightarrow+\infty}P_{\lambda,n}^{\omega}(\eta_{e^{C(\lambda)n}}\neq \emptyset)=1\text{\quad a.s.}
\end{equation}
with respect to $\mu$.
\end{theorem}

Theorem \ref{theorem 2.1 main} shows that our model has a critical value which is the inverse of ${\rm E}\rho{\rm E}\frac{1}{\xi}$. In the subcritical case, the process dies out before a moment with order $O(\log n)$ with high probability while in the supercritical case, the process survives at a moment with order $\exp(O(n))$ with high probability.

According to an intuitive idea, it is natural to guess that the critical value of the contact process is proportional to the mean of the recovery rate $\xi$, which holds trivially when $\xi$ is constant. However, our result shows that in the random environment the critical value is not proportional to the mean of $\xi$ simply but inversely proportional to the mean of the inverse of $\xi$, which we think is not without interest. In Section \ref{section 3} we will show that our result is consistent with a precise mean-field analysis.

\section{Mean field analysis}\label{section 3}
In this section, we give the mean-field analysis for this process. In brief, the mean-field analysis allows us to study a simpler model rather than the original process $\eta_t$. We first average the edge weights. That is, let all of the edge weights be ${\rm E}\rho.$ Then, we take the average of the recovery rates. In detail, suppose that the recovery rate $\xi$ takes value on ${Y}=\{y_1,y_2,\cdots,y_k\}$, where $y_j\in[1,M]$ and $k$ is a positive integer. Let $q_j:=\mu(\xi=y_j)$ for $j\in\{1,2,\cdots,k\}.$ Suppose also that, among the $n$ vertices, there are exactly $q_jn$ vertices with recovery rate $y_j$.  Let $A_t(j)$ be the number of the infected vertices which have recovery rate $y_j$ at time $t$. Please note that, by the above, an infected vertex with recovery rate $y_i$ becomes healthy at rate $y_j$, and a healthy vertex will be infected by an infected vertex at rate $\frac{\lambda}{n}{\rm E}\rho.$ Hence, we have
\begin{align}\label{equ 3.1}
\frac{d}{dt}A_t(i)=-y_iA_t(i)+\frac{\lambda}{n}{\rm E}\rho\sum_{j=1}^k[q_in-A_t(i)]A_t(j).
\end{align}
Let $f_t(i)=A_t(i)/q_in$ be the proportion of the infected ones in the $q_in$ vertices of recovery rate $y_i$ at time $t$. Then,by \eqref{equ 3.1}
\begin{align}\label{equ 3.2}
\frac{d}{dt}f_t(i)
&=-y_if_t(i)+\frac{\lambda}{n}{\rm E}\rho\sum_{j=1}^k[1-p_t(i)]A_t(j)\notag\\
&=-y_if_t(i)+{\lambda}{\rm E}\rho\sum_{j=1}^k[1-f_t(i)]f_t(j)q_j.
\end{align}

Now, suppose $\eta_t$ survives in an exponential length time, therefore, we may hope that the process will be stable or metastable in such a long time, which means we hope \eqref{equ 3.2} has a stable solution: $f_t(\cdot)=f(\cdot)$ for some $f(\cdot)\not\equiv0$. On the other hand, if $\eta_t$ dies out quickly, there should not be such a stable proportion $f(\cdot)$. So, we wonder whether there is a function $f(\cdot)\not\equiv0$ satisfies
\begin{align}\label{equ 3.3}
y_if(i)={\lambda}{\rm E}\rho\sum_{j=1}^k[1-f(i)]f(j)q_j.
\end{align}
Now, let $\sum_{j=1}^{k}f(j)q_j=x,$ then, by \eqref{equ 3.3}
\begin{align}\label{equ 3.4}
f(i)=\frac{\lambda x {\rm E}\rho}{y_i+\lambda x {\rm E}\rho}.
\end{align}
and then, by \eqref{equ 3.4}
\begin{align}\label{equ 3.5}
1=\sum_{j=1}^{k}\frac{\lambda q_j{\rm E}\rho}{y_j+\lambda x {\rm E}\rho}.
\end{align}
Note we hope to get some function $f(\cdot)\not\equiv0$, hence by \eqref{equ 3.4}, we should study whether there is some $\lambda$ such that equation \eqref{equ 3.5} has a positive solution.
Let $h(x)={\rm E}\frac{\lambda{\rm E}\rho}{\xi+\lambda x{\rm E}\rho},~x\in(0,+\infty),$ then \eqref{equ 3.5} is equivalent to $1=h(x)$. Note $h(x)$ is strictly monotonic decreasing in $(0,+\infty)$. By the bounded convergence theorem,
\[\lim_{x\rightarrow+\infty}h(x)=0,~{\rm and} \lim_{x\rightarrow0^+}h(x)=\lambda{\rm E}\rho{\rm E}\frac{1}{\xi}.\]
Therefore, define the critical value: $\lambda_c={({\rm E}\rho {\rm E}\frac{1}{\xi})}^{-1}$. By the above, when $\lambda>\lambda_c$, there exists an $x>0$ which solves \eqref{equ 3.5} and then $\eta_t$ survives in an exponential length time. When $\lambda<\lambda_c$, \eqref{equ 3.5} has no positive solution and therefore $\eta_t$ dies out quickly. Additionally, when $\lambda>\lambda_c$, by the monotonicity of $h(x)$, equation $1=h(x)$ and hence equation \eqref{equ 3.5} equivalently has an unique solution in $(0,+\infty)$ and we denote this unique solution by $x^*(\lambda)$ or $x^*$ for simplification. Hence, replacing $x$ in the right side of \eqref{equ 3.4} by $x^*$, we get the stable proportion and denote it by $f^*(\cdot).$

Until now, we have been taking a less rigorous approach to the phase transition of the process by the mean-field analysis. We will give the rigorous proof of our main results in the next two sections.

\section{Subcritical case}\label{section 4}
In this section we give the proof of \eqref{equ 2.1 subcritical}. Our approach is inspired by the graphical method introduced by Harris in \cite{Har1978}. In case of getting lost in the details, readers can read the conclusion at the end of this section first to obtain the logical procedure of this section.

First we introduce the graphical representation of our model. For each $n\geq 1$, we consider the graph $C_n\times [0,+\infty)$. That is to say, there is a time axis on each vertex on $C_n$. For any $1\leq i\neq j<+\infty$, we assume that $\{Y_i(t):t\geq 0\}$ is a Poisson process with rate $\xi(i)$ while $\{U_{(i,j)}(t):t\geq 0\}$ is a Poisson process with rate $\frac{\lambda}{n}\rho(i,j)$. Please note that here we care about the order of $i$ and $j$, so $U_{(i,j)}\neq U_{(j,i)}$. We assume that all these Poisson processes are independent. Intuitively, event moments of $Y_i(\cdot)$ are those at which vertex $i$ becomes healthy while event moments of $U_{(i,j)}(\cdot)$ are those at which vertex $i$ infects vertex $j$.

For each $1\leq i\leq n$ and any event moment $s$ of $Y_i(\cdot)$, we write a `$\Delta$' at $(i,s)$. For each $1\leq i\neq j\leq n$ and any event moment $r$ of $U_{(i,j)}(\cdot)$, we write an arrow `$\rightarrow$' from $(i,r)$ to $(j,r)$. For any $t>0$ and $i,j\in C_n$, we say that there is an infection path from $(i,0)$ to $(j,t)$ when there exists $m+1$ vertices $i=i_0, i_1, i_2,\ldots,i_m=j$ on $C_n$ and $m+2$ moments $0=t_{-1}<t_0<t_1<\ldots<t_m=t$ for some non-negative integer $m$ such that both the following conditions hold.

(1) For each $0\leq l\leq m-1$, there is an arrow `$\rightarrow$' from $(i_l,t_l)$ to $(i_{l+1},t_l)$.

(2) For each $0\leq l\leq m$, there is no `$\Delta$' at $\{i_l\}\times [t_{l-1},t_l]$.

Please note that we do not require that all the $m+1$ vertices $i_0,i_1,\ldots,i_m$ are different from each other in the above definition.

\begin{figure}
\centering
\includegraphics[scale=0.25]{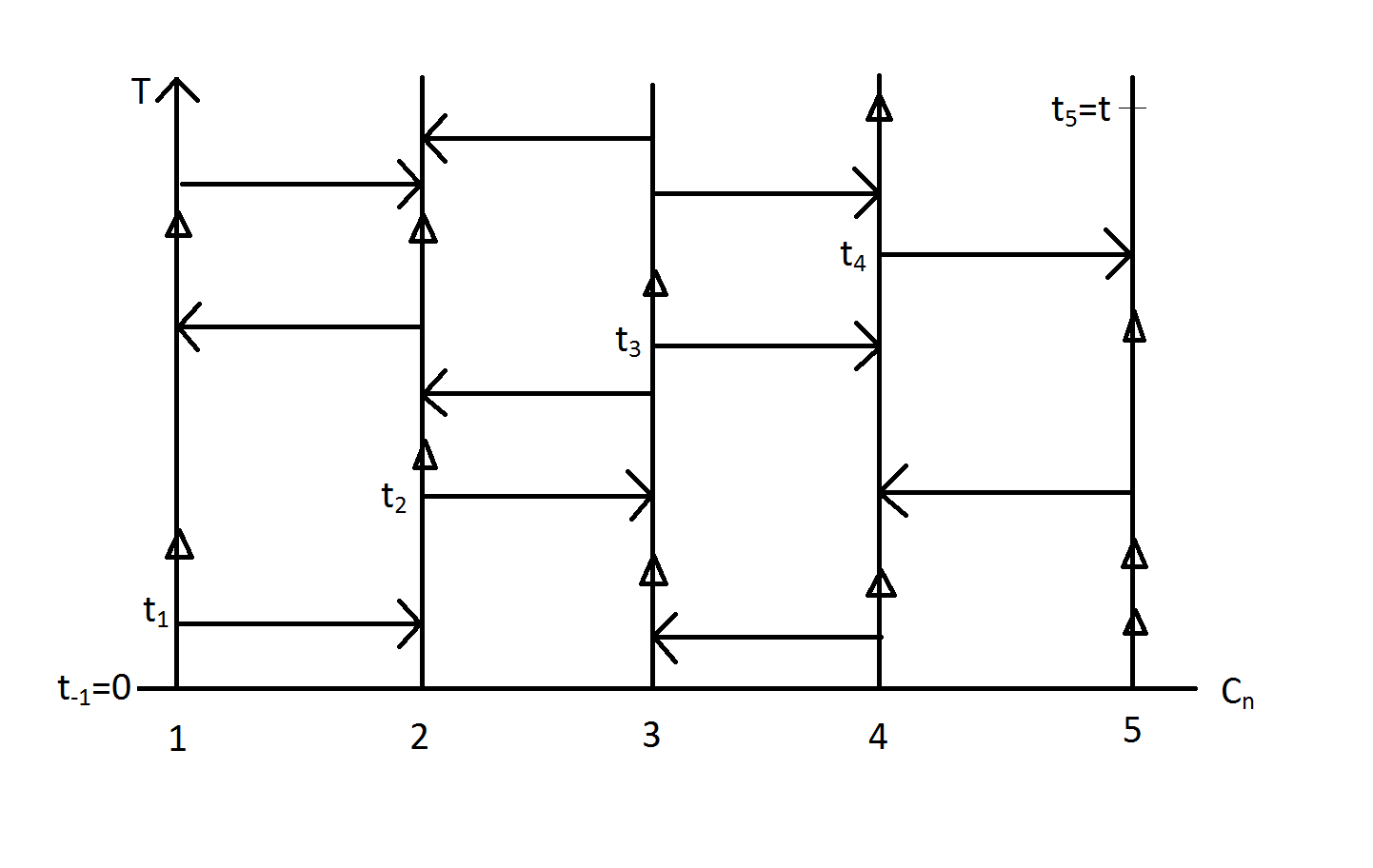}
\caption{Infection path}\label{figure 4.1}
\end{figure}

Figure \ref{figure 4.1} gives an example of an infection path. In Figure \ref{figure 4.1}, there is an infection path from $(1,0)$ to $(5,t)$.

According to the definition of our process and theory of graphical representation introduced in \cite{Har1978}, in the sense of coupling,
\begin{equation}\label{equ 4.1}
\{j\in \eta_t^A\}=\{\text{There is an infection path from $(i,0)$ to $(j,t)$ for some $i\in A$}\}.
\end{equation}
Therefore,
\[
\{j\in \eta_t^A\}=\bigcup_{i\in A}\{j\in \eta_t^{i}\}
\]
in the sense of coupling and hence
\begin{equation}\label{equ 4.2}
\{\eta_t^A\neq \emptyset\}=\bigcup_{i\in A}\{\eta_t^{i}\neq \emptyset\}.
\end{equation}
As a result,
\begin{equation}\label{equ 4.3}
P_{\lambda,n}(\eta_t^A\neq \emptyset)\leq |A| P_{\lambda,n}(\eta_t^{1}\neq \emptyset)
\end{equation}
according \eqref{equ 4.2} and the spatial homogeneity under the annealed measure. Please note that the symbol $\eta_t^{1}$ means that only vertex $1$ is infected at $t=0$ as we introduce in Section \ref{section 2}.

According to \eqref{equ 4.3} and Borel-Canteli Lemma, it is easy to prove that \eqref{equ 2.1 subcritical} is a direct corollary of the following lemma.
\begin{lemma}\label{lemma 4.1}
For any $\lambda<\frac{1}{{\rm E}\rho{\rm E}\frac{1}{\xi}}$, there exists $c(\lambda)>0$ such that
\begin{equation}\label{equ 4.4}
P_{\lambda,n}(\eta_{c(\lambda)\log n}^{1}\neq \emptyset)\leq
3n^{-3}.
\end{equation}
for sufficiently large $n$.
\end{lemma}

The main purpose of this section is to prove Lemma \ref{lemma 4.1}, but first we show how to utilize Lemma \ref{lemma 4.1} to prove \eqref{equ 2.1 subcritical}.

\proof[Proof of \eqref{equ 2.1 subcritical}]

According to Lemma \ref{lemma 4.1}, we choose $c(\lambda)$ satisfying \eqref{equ 4.4}, then by \eqref{equ 4.3} and \eqref{equ 4.4},
\begin{equation}\label{equ 4.5}
P_{\lambda,n}(\eta_{c(\lambda)\log n}\neq \emptyset)\leq
3nn^{-3}=3n^{-2}
\end{equation}
for sufficiently large $n$. For any $\epsilon>0$, by Chebyshev's inequality,
\begin{align}\label{equ 4.6}
\mu\big(\omega:P_{\lambda,n}^{\omega}(\eta_{c(\lambda)\log n}\neq \emptyset)>\epsilon\big)&\leq \frac{1}{\epsilon}{\rm E_\mu}\big(P_{\lambda,n}^{\omega}(\eta_{c(\lambda)\log n}\neq \emptyset)\big)\notag\\
&=\frac{1}{\epsilon}P_{\lambda,n}(\eta_{c(\lambda)\log n}\neq \emptyset).
\end{align}
By \eqref{equ 4.5} and \eqref{equ 4.6},
\begin{equation}\label{equ 4.7}
\sum_{n=1}^{+\infty}\mu\big(\omega:P_{\lambda,n}^{\omega}(\eta_{c(\lambda)\log n}\neq \emptyset)>\epsilon\big)<+\infty
\end{equation}
and hence
\begin{equation}\label{equ 4.8}
\mu\big(\omega:P_{\lambda,n}^{\omega}(\eta_{c(\lambda)\log n}\neq \emptyset)>\epsilon\text{~i.o.}\big)=0
\end{equation}
according to Borel-Cantelli Lemma. Therefore,
\begin{equation}\label{equ 4.9}
\limsup_{n\rightarrow+\infty} P_{\lambda,n}^{\omega}(\eta_{c(\lambda)\log n}\neq \emptyset)\leq \epsilon \text{~a.s.}
\end{equation}
with respect to $\mu$. Let $\epsilon\rightarrow 0$, then
\[
\lim_{n\rightarrow+\infty} P_{\lambda,n}^{\omega}(\eta_{c(\lambda)\log n}\neq \emptyset)=0 \text{~a.s.}
\]
with respect to $\mu$ and the proof is complete.

\qed

Now we only need to prove Lemma \ref{lemma 4.1}. The proof is divided into several steps. First we introduce the definition of an infection path with an given type and give an upper bound of $P_{\lambda,n}(\eta_t^{1}\neq\emptyset)$.

For any $m\geq 0$, we denote by $B_m$ the set of paths on $C_n$ starting at $1$ with length $m$. That is to say,
\begin{equation}\label{equ 4.10}
B_m=\{\vec{i}=(i_0,i_1,\ldots,i_m):i_0=1,1\leq i_l\leq n\text{~and~}i_{l-1}\neq i_l \text{~for~}1\leq l\leq m\}.
\end{equation}

For $t>0$, $\vec{i}=(i_0,i_1,\ldots,i_m)\in B_m$ and positive integers $j_0,j_1,j_2,\ldots,j_{m-1}$, we say that $\vec{i}$ is an infection path with type $(j_0,j_1,\ldots,j_{m-1})$ at moment $t$ when there exists $0=t_{-1}<t_0<t_1<\ldots<t_m=t$ such that all the following conditions hold.

(1) For each $0\leq l\leq m-1$, there is an arrow `$\rightarrow$' from $(i_l,t_l)$ to $(i_{l+1},t_l)$.

(2) For each $0\leq l\leq m$, there is no `$\Delta$' at $\{i_l\}\times [t_{l-1},t_l]$.

(3) For each $0\leq l\leq m-1$, $t_l$ is the `$j_l$'th event time after $t_{l-1}$ of $U_{(i_l,i_{l+1})}(\cdot)$.

\begin{figure}
\centering
\includegraphics[scale=0.25]{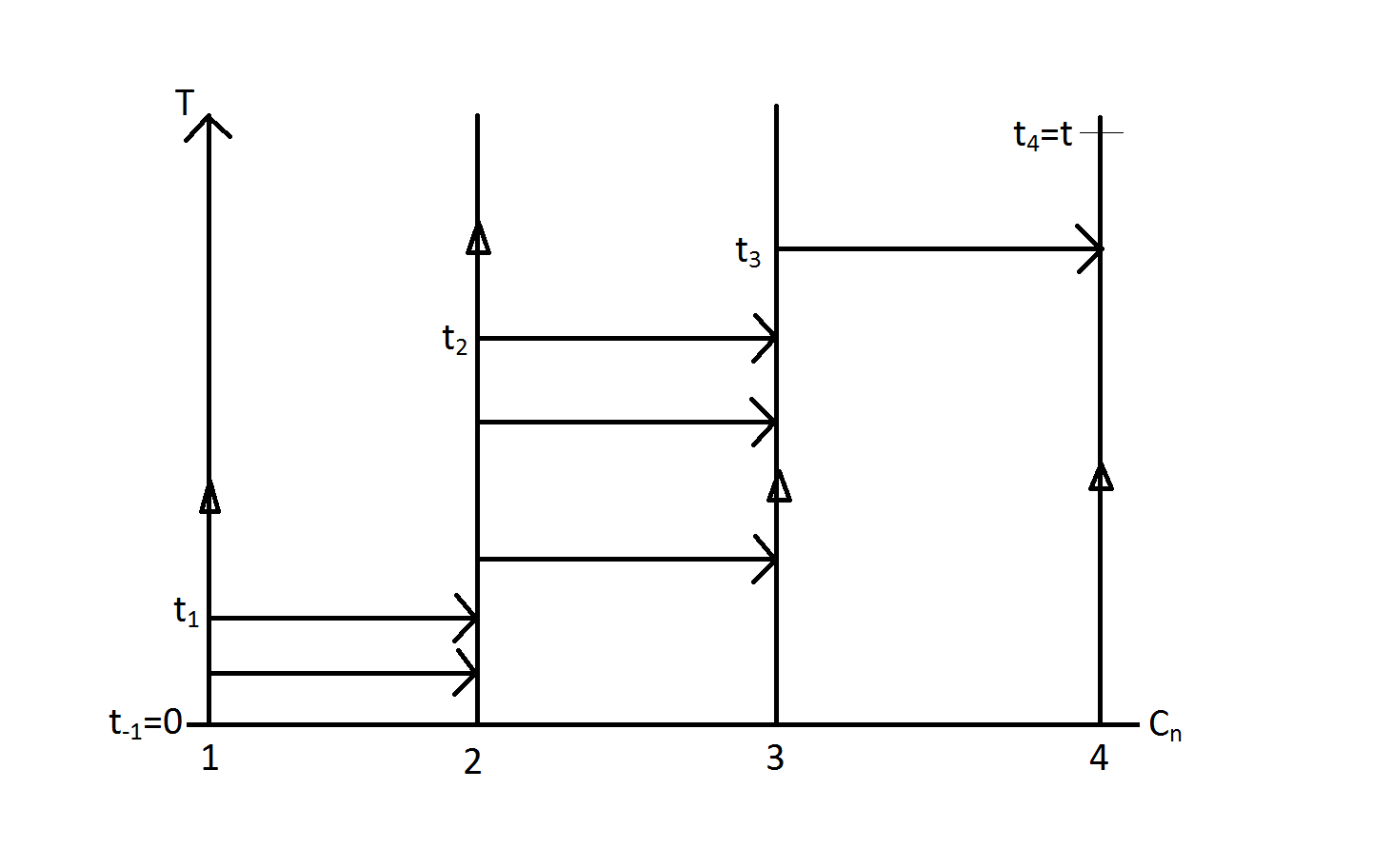}
\caption{infection path type}\label{figure 4.2}
\end{figure}

In Figure \ref{figure 4.2}, $(1,2,3,4)$ is an infection path with type $(2,3,1)$ at moment $t$. Please note that an infection path may be with more than one types. For example, in Figure \ref{figure 4.2}, $(1,2,3,4)$ is also with type $(1,3,1)$ at moment $t$.

For $\vec{i}\in B_m$, we denote by $A(\vec{i},j_0,j_1,\ldots,j_{m-1},t)$ the event that $\vec{i}$ is an infection path with type $(j_0,j_1,\ldots,j_{m-1})$ at moment $t$, then according to \eqref{equ 4.1},
\begin{equation}\label{equ 4.11}
\{\eta_t^{1}\neq \emptyset\}=\bigcup_{m=0}^{+\infty}\bigcup_{\vec{i}\in B_m}\bigcup_{j_0,j_1,\ldots,\atop j_{m-1}\geq 1}A(\vec{i},j_0,j_1,\ldots,j_{m-1},t).
\end{equation}
By \eqref{equ 4.11}, we have the following lemma which gives an upper bound of $P_{\lambda,n}^{\omega}(\eta_t^{1}\neq \emptyset)$.
\begin{lemma}\label{lemma 4.2}
For any $\omega\in \Omega$,
\begin{align}\label{equ 4.12}
&P_{\lambda,n}^{\omega}(\eta_t^{1}\neq \emptyset)
\leq \\
&\sum_{m=0}^{+\infty}(\frac{\lambda}{n})^{m}\sum_{\vec{i}\in B_m}\big(\prod_{l=0}^{m-1}\rho(\omega,i_l,i_{l+1})\big)\big(\prod_{l=0}^{m-1}\frac{1}{\xi(\omega,i_l)}\big)
P^{\omega,\vec{i}}(\sum_{l=0}^{m-1}V_l\leq t,\sum_{l=0}^mV_l\geq t),\notag
\end{align}
where under a probability measure $P^{\omega,\vec{i}}$ for each $\vec{i}=\{i_0,i_1,\ldots,i_m\}\in B_m$, $\{V_l\}_{l=0}^{m}$ are independent exponential times such that $V_l$ is with rate $\xi(\omega, i_l)$ for $0\leq l\leq m$ .
\end{lemma}

\proof

By \eqref{equ 4.11},
\begin{equation}\label{equ 4.13}
P_{\lambda,n}^{\omega}(\eta_t^{1}\neq \emptyset)\leq \sum_{m=0}^{+\infty}\sum_{\vec{i}\in B_m}\sum_{j_0,j_1\ldots,\atop j_{m-1}\geq 1}P_{\lambda,n}^{\omega}\big(A(\vec{i},j_0,j_1,\ldots,j_{m-1},t)\big).
\end{equation}
According to the definition of $A(\vec{i},j_0,j_1,\ldots,j_{m-1},t)$,
\begin{align}\label{equ 4.14}
&P_{\lambda,n}^{\omega}\big(A(\vec{i},j_0,j_1,\ldots,j_{m-1},t)\big)\notag\\
&=P^{\omega,\vec{i}}(\sum_{l=0}^{m-1}T_l\leq t, T_l\leq V_l \text{~for~}0\leq l\leq m-1, V_m\geq t-\sum_{l=1}^{m-1}T_l),
\end{align}
where $\{T_l\}_{l=0}^{m-1}$ are independent variables ,which are independent with $\{V_l\}_{l=1}^m$, such that $T_l$ is the sum of $j_l$ i.i.d exponential times with rate $\frac{\lambda}{n}\rho(\omega,i_l,i_{l+1})$. As a result, for $0\leq l\leq m-1$, $T_l$ has probability density function $p_l(t)$ given by
\[
p_l(t)=\frac{t^{j_l-1}\big(\frac{\lambda}{n}\rho(i_l,i_{l+1})\big)^{j_l}\exp\{-\big(\frac{\lambda}{n}\rho(i_l,i_{l+1})\big)t\}}{(j_l-1)!}.
\]
Therefore, by \eqref{equ 4.14},
\begin{align}\label{equ 4.15}
&P_{\lambda,n}^{\omega}\big(A(\vec{i},j_0,j_1,\ldots,j_{m-1},t)\big)\notag\\
&=\int\limits_{\sum\limits_{l=0}^{m-1}t_l\leq
t}\big(\prod_{l=0}^{m-1}p_l(t_l)\big)\big(\prod_{l=0}^{m-1}e^{-\xi(i_l)t_l}\big)e^{-\xi(i_m)(t-\sum\limits_{l=0}^{m-1}t_l)}dt_1dt_2\ldots
dt_{m-1}.
\end{align}
By \eqref{equ 4.13}, \eqref{equ 4.15} and repeatedly utilizations of
the fact that
\[
e^u=\sum_{l=0}^{+\infty}\frac{u^l}{l!},
\]
it is not difficult (but a little tedious) to check that the
right-hand side of \eqref{equ 4.13} equals
\begin{align}\label{equ 4.16}
&\sum_{m=0}^{+\infty}\big(\frac{\lambda}{n}\big)^m\sum_{\vec{i}\in
B_m}\Big(\prod_{l=0}^{m-1}\rho(i_l,i_{l+1})\Big)\notag\\
&\times
\int\limits_{\sum\limits_{l=0}^{m-1}t_l<t}\exp\big\{-\sum_{l=0}^{m-1}\xi(i_l)t_l-\xi(i_m)(t-\sum_{l=0}^{m-1}t_l)\big\}
dt_1dt_2\ldots dt_{m-1}.
\end{align}
According to the definition of $\{V_l\}_{l=1}^m$,
\begin{align}\label{equ 4.17}
&\int\limits_{\sum\limits_{l=0}^{m-1}t_l<t}\exp\big\{-\sum_{l=0}^{m-1}\xi(i_l)t_l-\xi(i_m)(t-\sum_{l=0}^{m-1}t_l)\big\}
dt_1dt_2\ldots dt_{m-1}\notag\\
&=\big(\prod_{l=0}^{m-1}\frac{1}{\xi(i_l)}\big)\int\limits_{\sum\limits_{l=0}^{m-1}t_l<t}(\prod_{l=0}^{m-1}\xi(i_l))\exp\big\{-\sum_{l=0}^{m-1}\xi(i_l)t_l\\
&-\xi(i_m)(t-\sum_{l=0}^{m-1}t_l)\big\}dt_1dt_2\ldots dt_{m-1}\notag\\
&=\big(\prod_{l=0}^{m-1}\frac{1}{\xi(i_l)}\big)P^{\omega,\vec{i}}(\sum_{l=0}^{m-1}V_l\leq
t,\sum_{l=0}^mV_l\geq t).\notag
\end{align}
Lemma \ref{lemma 4.2} follows from \eqref{equ 4.13}, \eqref{equ
4.16} and \eqref{equ 4.17} directly.

\qed

For $0\leq K\leq +\infty$, $\lambda,t>0$ and $n\geq 1$, we define
\begin{align*}
&F(\lambda,n,K,t)=\sum_{m=0}^{K}(\frac{\lambda}{n})^{m}\times\\
&\sum_{\vec{i}\in B_m}{\rm
E_\mu}\Big\{\big(\prod_{l=0}^{m-1}\rho(\omega,i_l,i_{l+1})\big)\big(\prod_{l=0}^{m-1}\frac{1}{\xi(\omega,i_l)}\big)
P^{\omega,\vec{i}}(\sum_{l=0}^{m-1}V_l\leq t,\sum_{l=0}^mV_l\geq
t)\Big\},
\end{align*}
then by Lemma \ref{lemma 4.2},
\begin{equation}\label{equ 4.18}
P_{\lambda,n}(\eta_t^{1}\neq \emptyset)\leq F(\lambda,n,+\infty,t).
\end{equation}

Our next goal is to show that Lemma \ref{lemma 4.1} follows from \eqref{equ 4.18} and the following three lemmas.

\begin{lemma}\label{lemma 4.3}
For any $\lambda>0$,let $\beta(\lambda)=\max\{2\lambda,1\}$. If
$\frac{C_2}{2}-C_1\log \beta(\lambda)>3$, then
\begin{equation}\label{equ 4.19}
F(\lambda,n,C_1\log n,C_2\log n)<n^{-3}
\end{equation}
for sufficiently large $n$.
\end{lemma}

\begin{lemma}\label{lemma 4.4}
For any $\lambda>0$ and $\theta>0$ such that
$\frac{M\lambda}{M+\theta}<\frac{1}{2}$, if $C_3\log 2-\theta
C_2>3$, then
\begin{equation}\label{equ 4.20}
F(\lambda,n,+\infty,C_2\log n)-F(\lambda,n,C_3\log n,C_2\log
n)<n^{-3}
\end{equation}
for sufficiently large $n$.

\end{lemma}

Please note that $M$ is the upper bound of $\xi$ as we introduced in
Section \ref{section 1}.

\begin{lemma}\label{lemma 4.5}
For any $\lambda<\lambda_c$, let
$\widehat{\lambda}=\frac{1+\lambda{\rm E}\rho{\rm
E}\frac{1}{\xi}}{2}$ and $C_3>C_1>3(\log
\frac{1}{\widehat{\lambda}})^{-1}$, then for sufficiently large $n$,
\begin{equation}\label{equ 4.21}
F(\lambda,n,C_3\log n,t)-F(\lambda,n,C_1\log n,t)<n^{-3}
\end{equation}
for any $t>0$.
\end{lemma}

We will give the proofs of these three lemmas later. Now we show how
to use these three lemmas to prove Lemma \ref{lemma 4.1}.

\proof[Proof of Lemma \ref{lemma 4.1}]

For $\lambda<\lambda_c$, we first choose $C_1(\lambda)$ satisfying
\[
C_1(\lambda)>3(\log \frac{1}{\widehat{\lambda}})^{-1},
\]
where
\[
\widehat{\lambda}=\frac{\lambda{\rm E}\rho{\rm
E}\frac{1}{\xi}+1}{2}.
\]

Then we choose $C_2(\lambda)$ such that
\[
\frac{C_2(\lambda)}{2}-C_1(\lambda)\log \beta(\lambda)>3,
\]
where
\[
\beta(\lambda)=\max\{2\lambda,1\}.
\]

Then we choose $C_3(\lambda)$ such that
\[
C_3(\lambda)\log 2-C_2(\lambda)\theta>3,
\]
where $\theta$ is sufficiently large such that
\[
\frac{M\lambda}{M+\theta}<\frac{1}{2}.
\]

According to Lemma \ref{lemma 4.5},
\begin{equation}\label{equ 4.22}
F(\lambda,n,C_3(\lambda)\log n,C_2(\lambda)\log
n)-F(\lambda,n,C_1(\lambda)\log n, C_2(\lambda)\log n)<n^{-3}
\end{equation}
for sufficiently large $n$.

According to Lemma \ref{lemma 4.4},
\begin{equation}\label{equ 4.23}
F(\lambda,n,+\infty,C_2(\lambda)\log n)-F(\lambda,n,C_3(\lambda)\log
n,C_2\log n)<n^{-3}
\end{equation}
for sufficiently large $n$.

According to Lemma \ref{lemma 4.3},
\begin{equation}\label{equ 4.24}
F(\lambda,n,C_1(\lambda)\log n,C_2(\lambda)\log n)<n^{-3}
\end{equation}
for sufficiently large $n$.

By \eqref{equ 4.22}, \eqref{equ 4.23} and \eqref{equ 4.24},
\begin{equation}\label{equ 4.25}
F(\lambda,n,+\infty,C_2(\lambda)\log n)<3n^{-3}
\end{equation}
for sufficiently large $n$.

Let $c(\lambda)=C_2(\lambda)$, then Lemma \ref{lemma 4.1} follows
from \eqref{equ 4.18} and \eqref{equ 4.25}.

\qed

Now we only need to prove Lemma \ref{lemma 4.3}, \ref{lemma 4.4} and
\ref{lemma 4.5}. First we prove Lemma \ref{lemma 4.3}.

\proof[Proof of Lemma \ref{lemma 4.3}]

Let $\{\widehat{V}_l\}_{l=0}^{+\infty}$ be i. i. d. exponential
times with rate $1$, then
\begin{equation}\label{equ 4.26}
P^{\omega,\vec{i}}(\sum_{l=0}^{m-1}V_l\leq t,\sum_{l=0}^mV_l\geq
t)\leq P(\sum_{l=0}^m\widehat{V}_l\geq t)
\end{equation}
since $\xi\geq 1$. By Chebyshev's inequality,
\begin{align}\label{equ 4.27}
&P(\sum_{l=0}^m\widehat{V}_l\geq C_2\log n)=P(\exp\{\frac{1}{2}\sum_{l=0}^m\widehat{V}_l\}\geq
n^{\frac{C_2}{2}})\notag\\
\leq & n^{-\frac{C_2}{2}}({\rm
E}e^{\frac{\widehat{V}_1}{2}})^{m+1}=n^{-\frac{C_2}{2}}2^{m+1}.
\end{align}

By \eqref{equ 4.26}, \eqref{equ 4.27} and the definition of
$F(\lambda,n,K,t)$,
\begin{align}\label{equ 4.28}
F(\lambda,n, C_1\log n, C_2\log n)\leq
2n^{-\frac{C_2}{2}}\sum_{m=0}^{C_1\log n}(2\lambda)^m,
\end{align}
since $\rho\leq 1,\xi\geq 1$ and
\[
|B_m|=(n-1)^m<n^m.
\]
Lemma \ref{lemma 4.3} follows from \eqref{equ 4.28} directly since
\[
\sum_{m=0}^{C_1\log n}(2\lambda)^m\leq K(\lambda)
\beta(\lambda)^{C_1\log n}=K(\lambda)n^{C_1\log\beta(\lambda)}
\]
for some constant $K(\lambda)>0$.

\qed

Now we give the proof of Lemma \ref{lemma 4.4}, which is similar
with that of Lemma \ref{lemma 4.3}.

\proof[Proof of Lemma \ref{lemma 4.4}]

Let $\{\widetilde{V}_l\}_{l=0}^{+\infty}$ be i. i. d. exponential
times with rate $M$, then
\begin{equation}\label{equ 4.29}
P^{\omega,\vec{i}}(\sum_{l=0}^{m-1}V_l\leq t,\sum_{l=0}^mV_l\geq
t)\leq P(\sum_{l=0}^{m-1}\widetilde{V}_l\leq t)
\end{equation}
since $\xi\leq M$.

By Chebyshev's inequality,
\begin{align}\label{equ 4.30}
&P(\sum_{l=0}^{m-1}\widetilde{V}_l\leq C_2\log n)=P(\exp\{-\theta\sum_{l=0}^{m-1}\widetilde{V}_l\}\geq
n^{-C_2\theta})\notag\\
&\leq n^{C_2\theta}({\rm E}e^{-\theta
\widetilde{V}_1})^m=n^{C_2\theta}(\frac{M}{M+\theta})^m
\end{align}
for any $\theta>0$.

By \eqref{equ 4.29}, \eqref{equ 4.30} and the definition of
$F(\lambda,n,K,t)$,
\begin{align}\label{equ 4.31}
&F(\lambda,n,+\infty,C_2\log n)-F(\lambda,n,C_3\log n,C_2\log
n)\notag\\
&\leq n^{C_2\theta}\sum_{m=C_3\log
n+1}^{+\infty}(\frac{M\lambda}{M+\theta})^m.
\end{align}
Please note that the facts that $\rho\leq 1,\xi\geq 1$ and
$|B_m|\leq n^m$ are also used to deduce \eqref{equ 4.31}.

For $\theta$ satisfying $\frac{M\lambda}{M+\theta}\leq \frac{1}{2}$,
\begin{equation}\label{equ 4.32}
\sum_{m=C_3\log n+1}^{+\infty}(\frac{M\lambda}{M+\theta})^m\leq
(\frac{1}{2})^{C_3\log n}=n^{-C_3\log 2}.
\end{equation}
Lemma \ref{lemma 4.4} follows from \eqref{equ 4.31} and \eqref{equ
4.32} directly.

\qed

To prove Lemma \ref{lemma 4.5}, we introduce simple random walk
$\{S_l\}_{l=0}^{+\infty}$ on $C_n$. In details, $S_0=1$ while
$S_{k}$ takes each vertex different with $S_{k-1}$ with probability
$\frac{1}{n-1}$ for each $k\geq 1$. For each $l\geq 1$, we denote by
$R_l$ the range of $\{S_0,S_1,\ldots,S_{l-1}\}$. That is to say,
\[
R_l=1+\sum_{k=1}^{l-1}1_{\{S_k\neq S_{k-1},S_k\neq S_{k-2},S_k\neq
S_{k-3},\ldots,S_k\neq S_0\}}.
\]
We denote by $\widehat{P}$ the probability measure of
$\{S_l\}_{l=0}^{+\infty}$ and $\widehat{{\rm E}}$ the expectation
operator with respect to $\widehat{P}$.

We introduce the following lemma to prove Lemma \ref{lemma 4.5}.
\begin{lemma}\label{lemma 4.6}
For any $\phi\in (0,1)$ and any $C_3>C_1>0$,
\begin{equation}\label{equ 4.33}
\lim_{n\rightarrow+\infty}\sup_{C_1\log n\leq m\leq C_3\log
n}\Big|\big(\widehat{\rm
E}\phi^{R_m}\big)^{\frac{1}{m}}-\phi\Big|=0.
\end{equation}
\end{lemma}

First we show how to use Lemma \ref{lemma 4.6} to prove Lemma
\ref{lemma 4.5}.

\proof[Proof of Lemma \ref{lemma 4.5}]

For $\lambda<\lambda_c$, we have $\lambda {\rm E}\rho{\rm
E}\frac{1}{\xi}<1$. Hence,
\begin{equation}\label{equ 4.34}
\widehat{\lambda}=\frac{\lambda {\rm E}\rho{\rm
E}\frac{1}{\xi}+1}{2}\in (\lambda {\rm E}\rho{\rm
E}\frac{1}{\xi},1).
\end{equation}
According to the definition of $\{S_l\}_{l=0}^{+\infty}$,
\begin{align}\label{equ 4.35}
&\sum_{\vec{i}\in B_m}{\rm
E_\mu}\Big\{\big(\prod_{l=0}^{m-1}\rho(\omega,i_l,i_{l+1})\big)\big(\prod_{l=0}^{m-1}\frac{1}{\xi(\omega,i_l)}\big)
P^{\omega,\vec{i}}(\sum_{l=0}^{m-1}V_l\leq t,\sum_{l=0}^mV_l\geq
t)\Big\}\notag\\
&\leq\sum_{\vec{i}\in B_m}{\rm
E_\mu}\Big\{\big(\prod_{l=0}^{m-1}\rho(\omega,i_l,i_{l+1})\big)\big(\prod_{l=0}^{m-1}\frac{1}{\xi(\omega,i_l)}\big)\Big\}\\
&=(n-1)^m\sum_{\vec{i}\in B_m}\frac{{\rm
E_\mu}\Big\{\big(\prod_{l=0}^{m-1}\rho(\omega,i_l,i_{l+1})\big)\big(\prod_{l=0}^{m-1}\frac{1}{\xi(\omega,i_l)}\big)\Big\}}{(n-1)^m}\notag\\
&=(n-1)^m\widehat{\rm E}\times {\rm
E}_\mu\Big\{\big(\prod_{l=0}^{m-1}\rho(\omega,S_l,S_{l+1})\big)\big(\prod_{l=0}^{m-1}\frac{1}{\xi(\omega,S_l)}\big)\Big\}.\notag
\end{align}
For each $m\geq 1$, there are $R_m$ different
vertices and at least $R_m-1$ different edges on the path $\{S_0,S_1,\ldots,S_{m-1}\}$. Therefore, according to
the fact that $\rho\leq 1$, $\xi\geq 1$ and the i. i. d. assumption
of the recovery rates and edge weights,
\begin{align}\label{equ 4.36}
{\rm
E}_\mu\Big\{\big(\prod_{l=0}^{m-1}\rho(\omega,S_l,S_{l+1})\big)\big(\prod_{l=0}^{m-1}\frac{1}{\xi(\omega,S_l)}\big)\Big\}
\leq \big({\rm E}\rho\big)^{R_m-1}\big({\rm
E}\frac{1}{\xi}\big)^{R_m}.
\end{align}
Let $C_1>3(\log\frac{1}{\widehat{\lambda}})^{-1}$ and $C_3>C_1$,
then by \eqref{equ 4.35}, \eqref{equ 4.36} and the definition of
$F(\lambda,n,K,t)$,
\begin{equation}\label{equ 4.37}
F(\lambda,n,C_3\log n,t)-F(\lambda,n,C_1\log n,t)\leq \frac{1}{{\rm
E}\rho}\sum_{m=C_1\log n+1}^{C_3\log n}\lambda^m\widehat{\rm E}({\rm
E}\rho{\rm E}\frac{1}{\xi})^{R_m}.
\end{equation}
By \eqref{equ 4.34}, we can choose $\delta$ sufficiently small such
that
\[
\lambda ({\rm E}\rho{\rm E}\frac{1}{\xi}+\delta)<\widehat{\lambda}.
\]
By Lemma \ref{lemma 4.6}, for sufficiently large $n$ and any
$C_1\log n\leq m\leq C_3\log n$,
\[
\lambda^m\widehat{\rm E}({\rm E}\rho{\rm E}\frac{1}{\xi})^{R_m} \leq
\lambda^m({\rm E}\rho{\rm
E}\frac{1}{\xi}+\delta)^m<\widehat{\lambda}^m.
\]

Hence by \eqref{equ 4.37}, for sufficiently large $n$,
\begin{align}\label{equ 4.38}
&F(\lambda,n,C_3\log n,t)-F(\lambda,n,C_1\log n,t)\leq \frac{1}{{\rm
E}\rho}\sum_{m=C_1\log
n}^{+\infty}\widehat{\lambda}^m \\
&\leq K_2(\lambda)\widehat{\lambda}^{C_1\log
n}=K_2(\lambda)n^{C_1\log \widehat{\lambda}}\notag
\end{align}
for some constant $K_2(\lambda)>0$, since $\widehat{\lambda}<1$.

Lemma \ref{lemma 4.5} follows from \eqref{equ 4.38} directly.

\qed

At last we give the proof of Lemma \ref{lemma 4.6}.

\proof[Proof of Lemma \ref{lemma 4.6}]

We write $\widehat{P}$ and $\widehat{\rm E}$ as $P$ and ${\rm E}$ in this proof since there is no misunderstanding. For any $j\geq 1$, we define
\[
\tau_j=\inf\{m\geq 0: R_m=j\},
\]
which is the first time when $\{S_l\}_{l=0}^{+\infty}$ visits at least $j$ different vertices. Then for any $\epsilon>0$,
\begin{equation}\label{equ 4.39}
P(R_m\leq m(1-\epsilon))=P(\tau_{(1-\epsilon)m}\geq m).
\end{equation}

Compared with the coupon collection model (See Section 2.5 of \cite{Dur2010}), it is not difficult to check that
\begin{equation}\label{equ 4.40}
\tau_{(1-\epsilon)m}\leq \sum_{l=1}^{(1-\epsilon)m}W_l
\end{equation}
in the sense of coupling, where $\{W_l\}_{l=1}^{+\infty}$ are i. i. d. random variables with geometric distribution
\[
P(W_1=k)=(\frac{m(1-\epsilon)}{n-1})^{k-1}(1-\frac{m(1-\epsilon)}{n-1}),
\]
for integer $k\geq 1$.

By \eqref{equ 4.39}, \eqref{equ 4.40} and Chebyshev's inequality,
\begin{equation}\label{equ 4.41}
P(R_m\leq m(1-\epsilon))\leq e^{-\theta m}\big({\rm E}e^{\theta W_1}\big)^{m(1-\epsilon)}
\end{equation}
for any $\theta>0$.

By \eqref{equ 4.41} and the distribution of $W_1$,
\begin{equation}\label{equ 4.42}
P(R_m\leq m(1-\epsilon))\leq \frac{e^{-\theta m\epsilon}}{\big\{1-\frac{me^\theta(1-\epsilon)}{n-1}\big\}^{m(1-\epsilon)}}.
\end{equation}

By \eqref{equ 4.42}, for any $\phi\in (0,1)$,
\begin{equation}\label{equ 4.43}
\phi^m\leq {\rm E}\phi^{R_m}\leq \phi^{(1-\epsilon)m}+\big[\frac{e^{-\theta \epsilon}}{\{1-\frac{m(1-\epsilon)e^\theta}{n-1}\}^{(1-\epsilon)}}\big]^m.
\end{equation}

We choose $\theta$ sufficiently large such that
\begin{equation}\label{equ 4.43 two}
2e^{-\theta\epsilon}<\phi^{1-\epsilon}.
\end{equation}
Since $\lim\limits_{n\rightarrow+\infty}\sup\limits_{C_1\log n\leq m\leq C_3\log n}|\{1-\frac{m(1-\epsilon)e^\theta}{n-1}\}^{(1-\epsilon)}-1|=0$,
\begin{equation}\label{equ 4.44}
\frac{1}{\{1-\frac{m(1-\epsilon)e^\theta}{n-1}\}^{(1-\epsilon)}}\leq 2
\end{equation}
for sufficiently large $n$ and each $C_1\log n \leq m\leq C_3\log n$.

By \eqref{equ 4.43}, \eqref{equ 4.43 two} and \eqref{equ 4.44},
\begin{equation}\label{equ 4.45}
\limsup_{n\rightarrow+\infty}\sup_{C_1\log n\leq m\leq C_3\log n}|({\rm E}\phi^{R_m})^{\frac{1}{m}}-\phi|\leq \phi^{1-\epsilon}-\phi.
\end{equation}

Let $\epsilon\rightarrow 0$, then Lemma \ref{lemma 4.6} follows from \eqref{equ 4.45}.

\qed

In conclusion, in this section we give the proof of \eqref{equ 2.1 subcritical}. First we show that \eqref{equ 2.1 subcritical} is a direct corollary of Lemma \ref{lemma 4.1}. To prove Lemma \ref{lemma 4.1}, we give Lemma \ref{lemma 4.2}-\ref{lemma 4.5} and show that Lemma \ref{lemma 4.1} follows from these four lemmas. We give the proofs of Lemma \ref{lemma 4.2}-\ref{lemma 4.5}, while the proof of Lemma \ref{lemma 4.5} is based on Lemma \ref{lemma 4.6}. The proof of Lemma \ref{lemma 4.6} is given at the end of this section.

\section{Supercritical case}\label{section 5}
In this section, we give the proof of \eqref{equ 2.2 supercritical}. The strategy of the proof comes from Peterson in \cite{Pet2011}. In \cite{Pet2011}, Peterson first handles with the case where the distribution of the recovery rate $\xi$ has finite support. Then, as a conclusion, get the proof of the general case. We will follow the same general outline to show \eqref{equ 2.2 supercritical} and what we do is modifying some details and some more sophisticated estimates.

Assume for now that the distribution of the recovery rate $\xi$ has finite support. We first give some notations which may have appeared in section \ref{section 3}. The recovery rates $\{\xi(i)\}_{i=1}^{+\infty}$ are random variables on probability space $\{\Omega,\mathfrak{F},\mu\}$ and in some finite space which is denoted by ${Y}=\{y_1,y_2,\cdots,y_k\}$, where $y_i\in[1,M]$ and and $k$ is a positive integer. Let $q_i:=\mu(\xi=y_i).$ As defined in section \ref{section 1}, $\forall \omega\in\Omega,$ $\eta_t$ is a contact process on $C_n$ that evolves in the ways of \eqref{equ 1.1 generator}. Now, given $j\in\{1,2,\cdots,k\}$, let
\begin{align}\label{equ 5.1}
A_t(j)=\sum_{m=1}^nI_{\{\xi(m)=y_j,\eta_t(m)=1\}}.
\end{align}
Hence, $A_t(j)$ is the number of the infected vertices which have recovery rate $y_j$ at time $t$. Let $A_t:=(A_t(1),A_t(2)\cdots,A_t(k))$. So, we want to show that the $k$ dimensional process $A_t$ will not visit the original point in an exponential length time whenever $\lambda>\lambda_c$.
Note $A_t$ takes value on $\mathbb{Z}_+^k$ and only one coordinate changes in each transition of $A_t$. When $\eta_t=\eta,$ let $A_t=A=(A(1),A(2),\cdots,A(k))\in \mathbb{Z}_+^k$ according to \eqref{equ 5.1}. Then, for $j\in\{1,2,\cdots,k\}$, $A\rightarrow(A(1),A(2),\cdots,A(j)+1,\cdots,A(k))$ at rate $q_j^+(\eta)$, and\\
$A\rightarrow(A(1),A(2),\cdots,A(j)-1,\cdots,A(k))$ at rate $q_j^-(\eta),$ where
\begin{equation}\label{equ 5.2}
\begin{cases}
&q_j^+(\eta)=\frac{\lambda}{n}\sum_{m=1}^n\big(I_{\{\xi(m)=y_j,\eta(m)=0\}}\sum_{i\not=m}\rho(i,m)I_{\{\eta(i)=1\}}\big),\\
&q_j^-(\eta)=y_j\sum_{m=1}^nI_{\{\xi(m)=y_j,\eta(m)=1\}}.
\end{cases}
\end{equation}

The following Lemma helps us bound the process $A_t$ from below. Before giving the Lemma, we first introduce some notations. Let
$S_n(i)=\sum_{m=1}^nI_{\{\xi(m)=y_i\}}$ be the number of vertices with recovery rate $y_i$ and $\hat{q}_n(i):=S_n(i)/n$ be the corresponding proportion. Please recall the proportion function $f(\cdot)$ defined in \eqref{equ 3.4} and the stable proportion $f^*(\cdot)$ in the end of section \ref{section 3}, then we introduce two associated $k$ dimensional sets. For $0<a< b< x^*,$ define
\[B_n(a,b):=\Big\{A\in\mathbb{Z}^k:\frac{\lambda a{\rm E}\rho}{y_i+\lambda a{\rm E}\rho}
\leq\frac{A(i)}{q_in}\leq\frac{\lambda b{\rm E}\rho}{y_i+\lambda b{\rm E}\rho},~i=1,2,\cdots,k\Big\},\] and
\[\widehat{B}_n(a,b):=\Big\{A\in\mathbb{Z}^k:\frac{\lambda a{\rm E}\rho}{y_i+\lambda a{\rm E}\rho}
\leq\frac{A(i)}{S_n(i)}\leq\frac{\lambda b{\rm E}\rho}{y_i+\lambda b{\rm E}\rho},~i=1,2,\cdots,k\Big\}.\]
Recalling that, in section \ref{section 3}, the stable proportion of infected vertices with recovery rate $y_i$ is $\frac{\lambda x^* {\rm E}\rho}{y_i+\lambda x^* {\rm E}\rho}.$ Hence, for $0<a< b< x^*,$ when $A_t\in B_n(a,b)$, hopefully, $A_t$ will have a positive drift in all of its coordinates with high probability. For $\eta \in\{0,1\}^{C_n},$ let $A(\eta)=(A(\eta,1),A(\eta,2)\cdots,A(\eta,k))$, where $A(\eta,i)=\sum_{m=1}^nI_{\{\xi(m)=y_i,\eta(m)=1\}}$.
Now, we give the following Lemma.
\begin{lemma}\label{lemma 5.1}
Suppose $\lambda>\lambda_c$. Then, $\forall a\in(0,x^*(\lambda)),$ there exists  $b\in(a,x^*(\lambda))$ and positive constants $\{\alpha_i^+,\alpha_i^-\}_{i=1}^k$ such that $\mu-a.s.,$
\[\sup_{\eta:A(\eta)=A}q_i^-(\eta)\leq\alpha_i^-n<\alpha_i^+n\leq \inf_{\eta:A(\eta)=A}q_i^+(\eta),~\forall A\in B_n(a,b),~i=1,2,\cdots,k.\] for all $n$ sufficiently large.
\end{lemma}
\proof[Proof of Lemma \ref{lemma 5.1}] We consider $\widehat{B}_n(a,b)$ first. By the definition of $\widehat{B}_n(a,b)$, we can obtain, when $A_t=A\in \widehat{B}_n(a,b),$ the number of all of the infected vertices is at least \[[\sum_{j=1}^k\frac{\lambda a{\rm E}\rho}{y_j+\lambda a{\rm E}\rho}S_n(j)],\] and the number of vertex with recovery rate $y_i$ which is not infected is at least $[\frac{y_i}{y_i+\lambda b{\rm E}\rho}S_n(i)].$ Hence, let
\begin{align}\label{equ 5.3}
b_i(A):=\frac{\lambda}{n}\inf\big\{\sum_{i\in E,j\in F}\rho({i,j}):
&|E|=[\sum_{j=1}^k\frac{\lambda a{\rm E}\rho}{y_j+\lambda a{\rm E}\rho}S_n(j)],\notag\\
&|F|=[\frac{y_i}{y_i+\lambda b{\rm E}\rho}S_n(i)],E\bigcap F=\emptyset\big\}
\end{align}
Then, by \eqref{equ 5.2} \[b_i(A)\leq \inf_{\eta:A(\eta)=A}q_i^+(\eta).\]
Note when $A(\eta)=A,$ by \eqref{equ 5.2}, we have: $q_i^-(\eta)=y_iA(i).$
Hence, let $d_i(A)=y_iA(i),$ then\[d_i(A)=\sup_{\eta:A(\eta)=A}q_i^-(\eta).\]

Let $A\in \widehat{B}_n(a,b),$ then,
\[d_i(A)\leq y_i\frac{\lambda b{\rm E}\rho}{y_i+\lambda b{\rm E}\rho}S_n(i)
=y_i\frac{\lambda b{\rm E}\rho}{y_i+\lambda b{\rm E}\rho}\hat{q}_n(i)n
:=\theta_i^-(a,b,n)n\]
By the law of large numbers, $\lim_{n\rightarrow\infty}\hat{q}_{n}(i)=q_i$. Hence, we have, $\mu-a.s.$,
\begin{align}\label{equ 5.4}
\lim_{n\rightarrow\infty}\theta_i^-(a,b,n)=y_i\frac{\lambda b{\rm E}\rho}{y_i+\lambda b{\rm E}\rho}q_i:=\theta_i^-(a,b)
\end{align}
Let $\theta_i^+(a,b,n):=b_i(A)/n$ and $\theta_i^+(a,b):=\lambda\cdot q_i\cdot \frac{y_i}{y_i+\lambda b{\rm E}\rho}\cdot {\rm E}\frac{\lambda a{\rm E}\rho}{\xi+\lambda a{\rm E}\rho}\cdot {\rm E}\rho$. Then, we claim that
\begin{align}\label{equ 5.5}
\lim_{n\rightarrow\infty}\theta_i^+(a,b,n)=\theta_i^+(a,b).
\end{align}
Let us continue with the proof of Lemma \ref{lemma 5.1} and the proof of \eqref{equ 5.5} will be in later. Recalling the function $h(x)$ defined in section \ref{section 3}, we have
\[\theta_i^+(a,b)/\theta_i^-(a,b)=\frac{a}{b}{\rm E}\frac{\lambda{\rm E}\rho}{\xi+\lambda a {\rm E}\rho}=\frac{a}{b}h(a).\]
Note, we suppose $\lambda>\lambda_c$. Hence, by section \ref{section 3}, equation $1=h(x)$ has an unique positive solution that is denoted by $x^*$. Since $h$ is decreasing and $h(x^*)=1,$ we have $h(a)>1$ since $a\in(0,x^*).$ Now , by taking proper $a'<a<b<b'<x^*(\lambda)$, we have \[\theta_i^+(a',b')/\theta_i^-(a',b')=\frac{a'}{b'}h(a')>1\].
Hence, we can take positive constants $\{\alpha_i^+,\alpha_i^-\}_{i=1}^k$ which satisfy
\[\theta_i^-(a',b')<\alpha_i^-<\alpha_i^+<\theta_i^+(a',b'),~{\rm for}~i=1,2,\cdots,k.\]
Therefore, by \eqref{equ 5.4}, \eqref{equ 5.5} and the law of large numbers,
\begin{align}\label{equ 5.6}
d_i(A)\leq\alpha_i^-n<\alpha_i^+n\leq b_i(A),~\forall A\in \widehat{B}_n(a',b'),~i=1,2,\cdots,k.
\end{align}
for $n$ sufficiently large. Therefore, take $b\in(a,b')$, then, the law of large numbers shows that: $B_n(a,b)\subset\widehat{B}_n(a',b')~{\rm when}~n~\rm{sufficiently~large}.$ Hence, by \eqref{equ 5.6}, we accomplish the proof of Lemma \ref{lemma 5.1} except the proof of \eqref{equ 5.5}.

We now turn to the proof of \eqref{equ 5.5}. Denote $s={\rm E}\frac{\lambda a{\rm E}\rho}{\xi+\lambda a{\rm E}\rho}$ and $t=q_i\cdot \frac{y_i}{y_i+\lambda b{\rm E}\rho}$. Hence, by the definition, $\theta_i^+(a,b)=st{\lambda}{\rm E}\rho.$ We first show that as $n\rightarrow\infty$,
\begin{align}\label{equ 5.7}
g_n(\omega):
&=\frac{\lambda}{n^2}\inf\big\{\sum_{i\in E,j\in F}\rho({i,j}):|E|=[sn],|F|=[tn],E\bigcap F=\emptyset\big\}\notag\\
&{\rightarrow}st\lambda{\rm E}\rho,~\mu-a.s..
\end{align}
Fix any $\delta>0$,
\begin{align}\label{equ 5.8}
&\sum_{n=1}^{\infty}\mu(\omega:|g_n(\omega)-st\lambda{\rm E}\rho|>\delta)\notag\\
\leq&\sum_{n=1}^{\infty}\mu(\omega:g_n(\omega)-st\lambda{\rm E}\rho>\delta)
+\sum_{n=1}^{\infty}\mu(\omega:g_n(\omega)-st\lambda{\rm E}\rho<-\delta)\notag\\
=&\textrm{I}+\textrm{II}.
\end{align}
By standard large deviation estimates, there exists positive constants $C_1$ such that
\begin{equation}\label{equ 5.9}
\textrm{I}\leq\sum_{n=1}^{\infty}\mu(\omega:\frac{\lambda}{n^2}\sum_{i=1}^{[sn]}\sum_{j=[sn]+1}^{[sn]+[tn]}
\rho(i,j)-st\lambda{\rm E}\rho>\delta)\leq\sum_{n=1}^{\infty}e^{-C_1n^2}<\infty,
\end{equation}
since
\[
g_n(\omega)\leq \frac{\lambda}{n^2}\sum_{i=1}^{[sn]}\sum_{j=[sn]+1}^{[sn]+[tn]}\rho(i,j).
\]
By standard large deviation estimates and Stirling's approximation, there exists positive constants $C_2$ and $C_3$ such that
\begin{align}\label{equ 5.10}
\textrm{II}&\leq\sum_{n=1}^{\infty}\dbinom{n}{[sn]}\dbinom{n-[sn]}{[tn]}\mu(\omega:\frac{\lambda}{n^2}\sum_{i=1}^{[sn]}\sum_{j=[sn]+1}^{[sn]+[tn]}\rho(i,j)-st\lambda{\rm E}\rho<-\delta)\notag\\
&\leq\sum_{n=1}^{\infty}\dbinom{n}{[sn]}\dbinom{n}{[tn]}e^{-C_3n^2}\\
&\leq\sum_{n=1}^{\infty}(C_2)^ne^{-C_3n^2}<\infty,\notag
\end{align}
since $g_n(\omega)<st\lambda {\rm E}\rho-\delta$ when and only when at least one pair of $E$ and $F$ satisfies that
\[
\frac{\lambda}{n^2}\sum_{i\in E,j\in F}\rho({i,j})<st\lambda {\rm E}\rho-\delta.
\]

Hence, \eqref{equ 5.7} follows from \eqref{equ 5.8}, \eqref{equ 5.9}, \eqref{equ 5.10} and Borel-Cantelli Lemma. Since when $|E|=[\sum_{j=1}^k\frac{\lambda a{\rm E}\rho}{y_j+\lambda a{\rm E}\rho}S_n(j)],$
$|E|/n\rightarrow s$ and when $|F|=[\frac{y_i}{y_i+\lambda b{\rm E}\rho}S_n(i)],$ $|F|/n\rightarrow t$, we obtain by \eqref{equ 5.3} and the definition of $\theta_i^+(a,b,n)$,
\[\lim_{n\rightarrow\infty}|\theta_i^+(a,b,n)-g_n(\omega)|=0,~\mu-a.s.\] Hence, by \eqref{equ 5.7}, we get \eqref{equ 5.5}. Thus we accomplish the proof of Lemma \ref{lemma 5.1}.

\qed\\
Define $Z_t=(Z_t(1),Z_t(2),\cdots,Z_t(k))$, where $\{Z_t(i)\}_{i=1}^k$ are independent nearest-neighbor random walks, and $Z_t(i)\rightarrow Z_t(i)+1$ at rate $\alpha_i^+n$ and $Z_t(i)\rightarrow Z_t(i)-1$ at rate $\alpha_i^-n$ for $i=1,2,\cdots,k$. The law of $Z_t$
with initial position $Z_0=Z$ will be denoted by $\widetilde{P}^{n,Z}.$ For any $a>0,$ define
\[U_n(a):=\Big\{Z=(Z(1),\cdots,,Z(k))\in\mathbb{Z}^k:{Z(i)}\geq\frac{\lambda a{\rm E}\rho}{y_i+\lambda a{\rm E}\rho}{q_in},~i=1,2,\cdots,k\Big\}\]
and $g(a)=(g_1^{a},g_2^{a},\cdots,g_k^{a}),$ where $g_i^{a}:=\lceil\frac{\lambda a{\rm E}\rho}{y_i+\lambda a{\rm E}\rho}{q_in}\rceil$.
Hence $g(a)$ is the smallest point in $U_n(a)$ under the partial order of the Euclidean space. Please note that, $\forall i\in\{1,2,\cdots,k\}$, $Z_t(i)$ is a nearest-neighbor random walk, whose positive drift rates are larger than negative drift rates and jump rates are proportional to $n$. Hence, we can get the following large deviation estimates about the process $Z_t$(we omit its proof here):
\begin{lemma}\label{lemma 5.2}
For the process $Z_t$, there exists a constant $C''>0$ such that $\forall\delta>0$ and $\forall t\geq0,$
\[\widetilde{P}^{n,g(\delta)}\big(Z_t\not\in U_n(\delta)\big)\leq ke^{-C''tn}.\]
\end{lemma}
\qed

For $\eta\in\{0,1\}^{C_n},$ define $X(\eta)=(X(\eta,1),X(\eta,2),\cdots,X(\eta,k)),$ where \[X(\eta,i)=\sum_{m=1}^nI_{\{\xi(m)=y_i,\eta(m)=1\}}.\]
Hence, $X(\eta_t)=A_t$ by \eqref{equ 5.1}. For $0<a<b$, define
\[\widetilde{B}_n(a,b):=\Big\{\eta\in\{0,1\}^{C_n}:\frac{\lambda a{\rm E}\rho}{y_i+\lambda a{\rm E}\rho}
\leq\frac{X(\eta,i)}{q_in}\leq\frac{\lambda b{\rm E}\rho}{y_i+\lambda b{\rm E}\rho},~i=1,2,\cdots,k\Big\},\]
and
\[\widetilde{U}_n(a):=\Big\{\eta\in\{0,1\}^{C_n}:\frac{X(\eta,i)}{q_in}\geq\frac{\lambda a{\rm E}\rho}{y_i+\lambda a{\rm E}\rho},~i=1,2,\cdots,k\Big\}.\]
For $0<t_1<t_2$, let $\eta[t_1,t_2]:=\{\eta_t:t\in[t_1,t_2]\}$ and $X(\eta[t_1,t_2]):=\{X(\eta_t):t\in[t_1,t_2]\}.$ We now give a Lemma which is crucial to the proof of \eqref{equ 2.2 supercritical}.
\begin{lemma}\label{lemma 5.3}
Suppose $\lambda>\lambda_c$. Let $\xi\in Y=\{y_1,y_2,\cdots,y_k\}$, where $y_i\in[1,M]$ and and $k$ is a positive integer. Then, for $a\in(0,x^*(\lambda))$, there exists a positive constant $C$ such that $\mu-a.s.,$
\begin{align}\label{equ 5.11}
\lim_{n\rightarrow\infty}P_{\lambda,n}^{\omega}\Big(\eta[0,e^{Cn}]\subset\widetilde{U}_n(a)\Big)=1.
\end{align}
In particular, \eqref{equ 5.11} implies that, \eqref{equ 2.2 supercritical} holds when the distribution of the recovery rate $\xi$ has finite support.
\end{lemma}
\proof[Proof of Lemma \ref{lemma 5.3}]It suffices to show that there exist positive constants $\gamma,$ $C'$ and $\delta\in(a,x^*(\lambda))$, which satisfy that $\mu-a.s.$,
\begin{align}\label{equ 5.12}
\inf_{\xi\in \widetilde{U}_n(\delta)}P_{\lambda,n}^{\omega}\big(\eta^{\xi}_{\gamma}\in \widetilde{U}_n(\delta),\eta^{\xi}[0,\gamma]\subset \widetilde{U}_n(a)\big)\geq1-e^{-C'n},
~{\rm when~}n~{\rm sufficiently~large.}
\end{align}
To see this, we can split the interval $[0,e^{C'n/2}]$ into $\frac{e^{C'n/2}}{\gamma}$ same parts with the same length $\gamma$, then by the Markov property,
\begin{align*}
&\lim_{n\rightarrow\infty}P_{\lambda,n}^{\omega}\Big(\eta[0,e^{C'n/2}]\subset \widetilde{U}_n(a)\Big)\\
&\geq\lim_{n\rightarrow\infty}(1-e^{-C'n})^{e^{C'n/2}/\gamma}\\
&\geq\lim_{n\rightarrow\infty}1-\frac{e^{-C'n/2}}{\gamma}=1,~~\mu-a.s..
\end{align*}
Hence, by taking $C=C'/2$, \eqref{equ 5.12} implies \eqref{equ 5.11} and we now turn to the proof of \eqref{equ 5.12}.
By the monotonicity of the process,
\begin{align*}
&\inf_{\xi\in \widetilde{U}_n(\delta)}P_{\lambda,n}^{\omega}\big(\eta^{\xi}_{\gamma}\in \widetilde{U}_n(\delta),\eta^{\xi}[0,\gamma]\subset \widetilde{U}_n(a)\big)\\
&=\inf_{\xi:X(\xi)=g(\delta)}P_{\lambda,n}^{\omega}\big(\eta^{\xi}_{\gamma}\in \widetilde{U}_n(\delta),\eta^{\xi}[0,\gamma]\subset \widetilde{U}_n(a)\big)\\
&\geq\inf_{\xi:X(\xi)=g(\delta)}P_{\lambda,n}^{\omega}\big(\eta^{\xi}_{\gamma}\in \widetilde{U}_n(\delta),\eta^{\xi}[0,\gamma]\subset\widetilde{B}_n(a,b)\big),
\end{align*}
where $g(\delta)$ is the smallest point in $U_n(\delta)$ under the partial order of Euclidean space as we defined before. Note, by the defination
\begin{align}\label{equ 5.13}
&\{\eta_{\gamma}\in \widetilde{U}_n(\delta)\}=\{X(\eta_{\gamma})\in U_n(\delta)\},\notag\\
{\rm and}~&\{\eta[0,\gamma]\subset \widetilde{B}_n(a,b)\}=\{X(\eta[0,\gamma])\subset B_n(a,b)\}.
\end{align}
Since $A_t=X(\eta_t)$, then, by Lemma \ref{lemma 5.1}, when $X(\eta_t)\in B_n(a,b)$, $X(\eta_t)$ can be dominated by $Z_t$ from below. Therefore, for any $\xi$ such that $X(\xi)=g(\delta)\in B_n(a,b),$ we have,
\begin{align}\label{equ 5.14}
\widetilde{P}^{n,g(\delta)}\big(Z_{\gamma}\in U_n(\delta)\big)
\leq & P_{\lambda,n}^{\omega}\big(X(\eta_{\gamma}^\xi)\in U_n(\delta),X(\eta^{\xi}[0,\gamma])\subset B_n(a,b)\big)\notag\\
&+ P_{\lambda,n}^{\omega}\big(X(\eta^{\xi}[0,\gamma])\not\subset B_n(a,b)\big).
\end{align}
Please note that, we can take proper positive constant $C$ (depends on $\lambda$ and $M$, independent of $n$) such that the total transition rate of $X(\eta_t)$ at any site is smaller than $Cn$ uniformly.
While, there exists also some constant $c$ such that the distance from $g(\delta)$ to the outside of $B_n(a,b)$ is larger than $cn$ (depends on $a,\delta,b$ and $M$, independent of $n$). Hence, if we let $Poisson(\theta)$ be Poisson distribution with parameter $\theta$, then
\begin{align}\label{equ 5.15}
P_{\lambda,n}^{\omega}\big(X(\eta^{\xi}[0,\gamma])\not\subset B_n(a,b)\big)\leq P(Poisson(\gamma Cn)\geq cn).
\end{align}
 Note that $Poisson(\gamma Cn)$ can be seen as $n$ independent Poisson distributions with parameter $\gamma C$, then by the standard large deviation estimates, we can choose $\gamma\in(0,c/C)$ and $C'''>0,$ such that
\begin{align}\label{equ 5.16}
P(Poisson(\gamma Cn)\geq cn)\leq e^{-C'''n}.
\end{align}
Hence, by \eqref{equ 5.13}, and \eqref{equ 5.14},
\begin{align}\label{equ 5.17}
\widetilde{P}^{n,g(\delta)}\big(Z_{\gamma}\in U_n(\delta)\big)
\leq & P_{\lambda,n}^{\omega}\big(\eta^{\xi}_{\gamma}\in \widetilde{U}_n(\delta),\eta^{\xi}[0,\gamma]\subset\widetilde{B}_n(a,b)\big)\notag\\
&+ P_{\lambda,n}^{\omega}\big(\eta^{\xi}[0,\gamma]\not\subset \widetilde{B}_n(a,b)\big).
\end{align}
By \eqref{equ 5.13}, \eqref{equ 5.15}, and \eqref{equ 5.16},
\[P_{\lambda,n}^{\omega}\big(\eta^{\xi}[0,\gamma]\not\subset \widetilde{B}_n(a,b)\big)\leq e^{-C'''n}.\]
Thus, by Lemma \ref{lemma 5.2}, \eqref{equ 5.17} and taking $C'\in(0,\min\{C'',C'''\})$, we get that $\mu-a.s.$,
\[
P_{\lambda,n}^{\omega}\big(\eta^{\xi}_{\gamma}\in \widetilde{U}_n(\delta),\eta^{\xi}[0,\gamma]\subset\widetilde{B}_n(a,b)\big)
\geq1-e^{-C'''n}-ke^{-C''\gamma n}\geq1-e^{-C'n},
\]
for all $n$ large enough. Hence, we get the proof of \eqref{equ 5.12} and then accomplish the proof of Lemma \ref{lemma 5.3}.

\qed

We now turn to the case where the distribution of $\xi$ does not have finite support.
\proof[Proof of \eqref{equ 2.2 supercritical}, general case]We first introduce the finite approximations to the process. In detail, for each $m\geq1,$ we modify $\{\xi(i)\}_{i\geq1}$ by $\{\xi_m(i)\}_{i\geq1}$, where $\xi_m(i)=\frac{\lceil \xi(i)m\rceil}{m}$. For each $m\geq1,$ $\eta_t^{m}$ is a contact process on $C_n$ defined in the ways of \eqref{equ 1.1 generator} and \eqref{equ 1.2 spin function} but with random recovery rates $\{\xi_m(i)\}_{1\leq i\leq n}$ and edge weights $\{\rho(i,j)\}_{1\leq i\neq j\leq n}.$ Since $\xi(i)\in[1,M],$ the distribution of $\xi_m(i)$ has finite support. Hence, \eqref{equ 5.11} holds for the process $\eta_t^{m}$ by Lemma \ref{lemma 5.3}. Define $\lambda_c^{(m)}:=({\rm E}\rho{\rm E}\frac{1}{\xi_m(1)})^{-1}.$ Hence, $\lambda_c^{(m)}$ is the critical value for the process $\eta_t^{m}$. Since $\lim_{m\rightarrow\infty}\xi_m(1)=\xi(1),$ then, by the  bounded convergence theorem
\[\lim_{m\rightarrow\infty}\lambda_c^{(m)}=\lambda_c.\] Thus, when $\lambda>\lambda_c,$ we can take $m$ sufficiently large so that $\lambda>\lambda_c^{(m)}.$ Hence, for such fixed $m$, consider the process $\eta_t^m$. By Lemma \ref{lemma 5.3}, there exists a positive constant $C$ such that $\mu-a.s.$,
\begin{align}\label{equ 5.18}
P_{\lambda,n}^{\omega}(\eta^m_{e^{Cn}}\not=\emptyset)\rightarrow1,~{n\rightarrow\infty.}
\end{align}
Note the recovery rate $\xi_m(i)\geq\xi(i)$ for all $i\geq1$, then by the monotonicity of the contact process, $\eta_t$ is always larger than $\eta_t^m$ in the sense of coupling. Therefore, we have
\[P_{\lambda,n}^{\omega}(\eta_{e^{Cn}}\not=\emptyset)\geq P_{\lambda,n}^{\omega}(\eta^m_{e^{Cn}}\not=\emptyset).\]
Hence, by \eqref{equ 5.18}, we complete the proof of \eqref{equ 2.2 supercritical}.

\qed

\textbf{Acknowledgments.} The authors are grateful to the financial support from the National Natural Science Foundation of China with grant number 11371040, 11531001 and 11501542.

{}
\end{document}